\documentclass[12pt]{article}

\usepackage{amsmath}
\usepackage{lipsum}

\usepackage{amsthm}
\usepackage{mathrsfs}	
\usepackage{amssymb}
\usepackage{amscd}
\usepackage{setspace}	
\usepackage{tikz}		
\usepackage{tikz-cd}
\usetikzlibrary{matrix}	
\usepackage{graphicx}	
\graphicspath{ {pics/} }	
\usepackage{wrapfig}	
\usepackage{enumerate}	
\usepackage{bbm}		

\usepackage[utf8]{inputenc}	
\usepackage[english]{babel}	
\usepackage{blindtext}		
\usepackage[affil-it]{authblk}	

\usepackage{faktor}			
\usepackage{xparse}			

\usepackage{hyperref}		
\usepackage{cite}			

\newtheorem{theorem}{Theorem} 

\newtheorem{proposition}[theorem]{Proposition}

\newtheorem{corollary}[theorem]{Corollary}

\newcommand{\acts}{\curvearrowright}

\newcommand{\C}{\mathbb{C}}

\begin{document}
\title{Isometric Actions are Quasidiagonal}
\author{Samantha Pilgrim}
\maketitle

\paragraph{Abstract:}  We show every isometric action is quasidiagonal in a strong sense.  This shows that reduced crossed products by such actions are quasidiagonal or MF whenever the reduced algebra of the acting group is quasidiagonal or MF.  

\vspace{3mm}
Quasidiagonality of group actions was first introduced in \cite[3.2]{kerr2011residually}.  The proof of \ref{translation actions theorem} relies on convolution with an approximate identity having certain properties.  We begin with a construction to produce such kernels.  

\proposition{Suppose $G$ is a compact group.  Then there exists a net of functions $k_i\in C(G)$ such that the following hold for all $i$:
\begin{description}
\item[(1)] $k_i\geq 0$
\item[(2)] $\int_G k_i d\mu = 1$
\item[(3)] $k_i * f$ converges uniformly to $f$ for all $f\in C(G)$
\item[(4)] $k_i$ lies in a finite-dimensional $G$-invariant subspace of $C(G)$
\end{description}
}

\begin{proof}
First, consider a neighborhood basis $(U_i)$ for the identity in $G$ and let $(h_i)$ be a net of non-negative, $L^1(G)$-normalized continuous functions with $h_i$ supported inside $U_i$.  Then $(h_i)$ is an approximate identity in the sense that $f*h_i\to f$ uniformly for $f\in C(G)$.  So $(h_i)$ satisfies $(1)$, $(2)$, and $(3)$.  

Let $C_{\text{fin}}(G)$ be the collection of functions in $C(G)$ whose $G$-orbits are contained in a finite-dimensional subspace.  Then $C_{\text{fin}}(G)$ is dense in $C(G)$ by the Peter-Weyl theorem (as it contains the matrix coefficients of all irreducible unitary representations of $G$).  Also, $C_{\text{fin}}(G)$ is a sub $*$-algebra: if $f$ and $g$ are linear combinations of $\{f_i\}_{i=1}^n$ and $\{g_j\}_{j=1}^m$ respectively, then $fg$ is a linear combination of functions in $\{f_ig_j\}_{i, j}$, and the action $G\acts C(G)$ respects products so that a product of functions with finite orbits will also have a finite orbit.  

Now, let $\sqrt{h_i}$ be the pointwise square root of $h_i$.  Let $l_i\in C_{\text{fin}}(G)$ be such that $\| l_i - \sqrt{h_i}\|_{C(G)}\leq \frac{\epsilon_i}{\|\sqrt{h_i}\|}_{C(G)}$ with $\epsilon_i\to 0$.  Put $k_i = \frac{l_il_i^*}{\|l_il_i^*\|}_1$.  

Note that $k_i$ satisfies $(1)$ by construction.  Moreover, $\|l_il_i^* - h_i\|_{C(G)}\leq \frac{\epsilon_i}{\|\sqrt{h_i}\|_{C(G)}} \cdot 2\|\sqrt{h_i}\|_{C(G)} + \frac{\epsilon_i^2}{\|\sqrt{h_i}\|^2_{C(G)}} < 3\epsilon_i\to 0$ as $i\to \infty$,  so $\Big|\int l_il_i^*d\mu - 1\Big| < 3\epsilon_i$.  But that means normalizing $l_il_i^*$ to have integral equal to $1$ scales it by a real number between $\frac{1}{1 + 3\epsilon_i}$ and $\frac{1}{1 - 3\epsilon_i}$ with $\epsilon_i\to 0$.  Thus, $(k_i)$ satisfies $(3)$ since $\|l_il_i^* - h_i\|_{C(G)}\to 0$ and $\|k_i - l_il_i^*\|_{C(G)}\to 0$.  It also satisfies $(4)$ since each $l_i$ satisfies $(4)$ and $C_{\text{fin}}(G)$ is a $*$-subalgebra.  \end{proof}


\theorem{Suppose $G$ is a compact group.  Then the action $G\acts C(G)$ by left multiplication is quasidiagonal.  In fact, the u.c.p. maps can be taken to be approximately equivariant for all $\gamma\in G$.  }\label{translation actions theorem}
\begin{proof}
Fix $F\subset C(G)$ a finite subset and $\epsilon>0$.  We can assume each $f\in F$ has $\|f\|_{C(G)}\leq 1$.  Let $(k_i)$ be as in the previous proposition.  Let $W_i$ be a finite dimensional, $G$-invariant subspace of $C(G)$ containing $k_i$.  

Denote by $\Phi_i$ the map $C(G)\to C(G)$ given by $f\mapsto f*k_i$.  This is equivariant for the left translation action.  Also the image of $\Phi_i$ lies in $W_i$.  To see this, observe that $(f*k_i)(g) = \int_G f(h)k_i(h^{-1}g)dh = \int_G f(h)(\lambda_hk_i)(g)dh = \Big(\int_G f(h)(\lambda_hk_i)dh\Big)(g)$ where $\lambda$ is the left translation action and the last integral in parenthesis is $C(G)$-valued.  This integral is then a limit of linear combinations of $(\lambda_hk_i)$ for different $h\in G$, and so is contained in $W_i$ since $W_i$ is a closed, invariant subspace.  Moreover, since $k_i\geq 0$, $\Phi_i$ is positive, hence completely positive as a map $C(G)\to C(G)$. 


Consider $W_i$ now with the $L^2(G)$-inner product.  This is a finite dimensional Hilbert space on which $G$ acts unitarily.  Finite-dimensionality implies we can find a sufficiently large $E_i\subset G$ such that the map $\psi_i: W_i\to \C^{E_i}$ given by taking evaluations is an isomorphism of $W_i$ onto its image and preserves the suprema of functions in $(F\cup F^2)*k_i \pm [(F\cup F^2)*k_i]^2$ (and we assume $0\in F$), since continuous functions on $G$ achieve their suprema.  Then we can conjugate over the action on $W_i$ to an action on a subspace of $\C^{E_i}$ and then extend trivially to the orthogonal complement of this subspace to get an action (representation) on all of $\C^{E_i}$.  Since all of $C(G)$ is mapped by $\psi_i\circ \Phi_i$ into the subspace, extending this way won't affect our estimates later.  

Partition $G$ into measurable sets, $P_e$ each containing exactly one $e\in E_i$.  For a sufficient choice of $E_i$, we can ensure the diameters of the $P_e$ are as small as we wish.  Equip $\C^{E_i}$ with the inner product $\langle v, w\rangle_i = \sum_{e\in E_i}v(e)\overline{w(e)}\mu(P_e)$ where $\mu$ is normalized Haar measure on $G$.  Notice that the canonical basis of $\C^{E_i}$ is still orthogonal with respect to this inner product, so functions on $E_i$, thought of as multiplication operators, are still represented by the same diagonal matrices as they would be with respect to the canonical basis.  

Since $\psi_i$ is not unitary as a map $(C(G), \langle\cdot, \cdot\rangle_{L^2G})\to (\C^{E_i}, \langle\cdot, \cdot\rangle_i)$, the representation of $G$ on $(\C^{E_i}, \langle\cdot, \cdot\rangle_i)$ given above is not unitary.  However, since $G$ is compact, we can make this representation unitary by replacing the original inner product with its average, $\langle\cdot, \cdot\rangle_i^*$, over $G$ .  We now have a unitary representation of $G$ on $(\C^{E_i}, \langle\cdot, \cdot\rangle_i^*)$, hence an action of $G$ on $B(\C^{E_i}, \langle\cdot, \cdot\rangle_i^*)$.  We can also think of $\C^{E_i}$ as multiplication operators on $(\C^{E_i}, \langle\cdot, \cdot\rangle_i^*)$ and define $\tilde{\psi_i}: W_i\to B(\C^{E_i}, \langle\cdot, \cdot\rangle_i^*)$.  

Some care is needed here since, although the averaged inner product is equivalent to the original, a priori the constants could increase as $i$ increases so that our errors are amplified.  However, if we assume $E_i$ is sufficiently large and that diameters of the $P_e$ are sufficiently small, we have that $\langle f, g\rangle \approx_\epsilon \langle \psi_i(f), \psi_i(g)\rangle_i$ for all $f, g\in B_2(W_i)$.  In other words, $\psi_i$, and hence the action on $(\C^{E_i}, \langle\cdot, \cdot\rangle)$, are `almost unitary'.   Assume $\epsilon <1$ so that the pullback of the unit ball by $\psi_i$ is contained in $B_2(W_i)$.  Then for $\gamma\in G$ and $\psi_i(f),\psi_i(g)\in B_1(\psi_i(W_i))$, $\langle \gamma\cdot \psi_i(f), \gamma\cdot \psi_i(g)\rangle_i = \langle \psi_i(\gamma\cdot f), \psi_i(\gamma\cdot g)\rangle_i\approx_\epsilon \langle \gamma\cdot f, \gamma\cdot g\rangle_i = \langle f, g\rangle_i\approx_\epsilon \langle \psi_i(f), \psi_i(g)\rangle$.  This shows $\langle v, w\rangle_i^* \approx_{2\epsilon} \langle v, w\rangle_i$ for $v, w\in B_1(\C^{E_i})$ (since the two are the same on the orthogonal complement of $\psi_i(W_i)$) and so errors are only amplified by a constant multiple which tends to $1$ as $E_i$ becomes larger.  

A second issue arises since the canonical basis of $\C^{E_i}$ is no longer orthonormal with respect to $\langle \cdot, \cdot \rangle_i^*$.  We therefore have to choose a new basis, which means the multiplication operators coming from $\C^{E_i}$ are represented by different matrices.  This means the map $\tilde{\psi_i}: W_i\to B(\C^{E_i}, \langle\cdot, \cdot\rangle_i^*)$ may not be $*$-linear, since $\tilde{\psi_i}(f) = V\psi_i(f)V^{-1}$ where $V$ is invertible, but not unitary, and we think of $\psi_i(f)$ as a diagonal matrix.  Notice that, for an appropriate choice of $E_i$, $V$ takes a basis which is orthonormal (with respect to $\langle\cdot, \cdot\rangle_i^*$) to a basis which is orthonormal with respect to $\langle\cdot, \cdot\rangle$ and hence has $\langle e_j, e_k\rangle_i^* \approx_{2\epsilon} 0$ and $\langle e_j, e_j\rangle_i^*\approx_{2\epsilon} 1$.  This implies $V$ is `almost unitary' in the sense that $\langle Vv, Vu\rangle_i^*\approx_{4\epsilon} \langle v, u\rangle_i^*$.  Recall that we can write $V = AU$ in a unique way where $A$ is a positive matrix and $U$ is a unitary.  Since $V$ is invertible and bounded below and above by $1 - 4\epsilon$ and $1 + 4\epsilon$ respectively, the same is true for $A$.  Thus, $A$ is within $4\epsilon$ of the identity.  The matrix $U$ is therefore a unitary such that $\|V - U\|<4\epsilon$.  To simplify things later, we can go back and choose $E_i$ so that  $\|V - U\|$ is sufficiently small so that $\|V\psi_i(f)V^{-1} - U\psi_i(f)U^*\|<\epsilon$ for all $f$ with $\|f\|_{C(G)}\leq 1$.

Defining $\Psi_i(f) := U\psi_i(f*k_i)U^*$ now gives a map $C(G)\to B(\C^{E_i}, \langle\cdot, \cdot\rangle_i^*)$ which is a composition of a unital, completely positive map and a unital $*$-homomorphism and hence u.c.p..  Since $\Psi_i$ is an $\epsilon$-perturbation of an equivariant map, $\|\Psi_i(\gamma\cdot f) - \gamma\cdot \Psi_i(f)\|_{B(\C^{E_i}, \langle\cdot, \cdot\rangle^*)}<\epsilon$ for all $f\in F$ and $\gamma\in G$.  Notice that this approximation works for all of $G$ rather than just a finite subset.  

Moreover, for $f\in F$, we can choose $i$ and $E_i$ so that 

$$\|f\|_{C(G)}\approx_{\epsilon} \|\psi_i(f*k_i)\|_{C(G)} = \|\psi_i(f*k_i)\|_{B(\C^{E_i}, \langle\cdot, \cdot\rangle)} \approx_\epsilon \|V\psi_i(f*k_i)V^{-1} \|_{B(\C^{E_i}, \langle\cdot, \cdot\rangle^*)} $$

$$\approx_\epsilon \|U\psi_i(f*k_i)U^* \|_{B(\C^{E_i}, \langle\cdot, \cdot\rangle^*)} = \|\Psi_i(f*k_i)\|_{B(\C^{E_i}, \langle\cdot, \cdot\rangle^*)}$$

%

for all $f\in F$.  


Also, if we choose $i$ large enough that $f, g,$ and $fg$ are all approximated to within $\epsilon$ by their convolutions with $k_i$, we have

%

$$\|U\psi_i((fg)*k_i)U^* - U\psi_i(f*k_i)U^*U\psi_i(f*k_i)U^*\|_{B(\C^{E_i}, \langle\cdot, \cdot\rangle_i^*)}$$
$$= \|U\psi_i((fg)*k_i) - \psi_i(f*k_i)\psi_i(g*k_i)U^*\|_{B(\C^{E_i}, \langle\cdot, \cdot\rangle_i^*)} = \|U\psi_i((fg)*k_i - (f*k_i)(g*k_i))U^*\|_{B(\C^{E_i}, \langle\cdot, \cdot\rangle_i^*)}$$
$$=\|(fg)*k_i - (f*k_i)(g*k_i)\|_{C(G)} = \|\Big(((fg)*k_i - fg) + fg\Big) - ((f*k_i - f) + f)((g*k_i - g) + g)\|_{C(G)}$$
$$\leq \epsilon + \|fg - ((f*k_i - f)(g*k_i - g) + (f*k_i - f)(g) + (g*k_i - g)f + fg)\|_{C(G)}$$
$$\leq \epsilon + \epsilon^2 + \epsilon + \epsilon < 4\epsilon$$
using that conjugation by $U$ and $\psi_i$ are both multiplicative and that the moduli of $f$ and $g$ are bounded above by $1$.  




Observing that $B(\C^{E_i}, \langle\cdot, \cdot\rangle_i^*)$ identifies with $M_{|E_i|}(\C)$ now completes the proof.  \end{proof}

\normalfont

It is now more or less straightforward to extend the previous theorem and obtain our main result.  

\theorem{Suppose $\Gamma\acts X$ is an isometric action by a countable discrete group on a compact space.  Then $\Gamma\acts X$ is quasidiagonal.  In fact, the u.c.p. maps can be taken to be approximately equivariant for all $\gamma\in \Gamma$.  }\label{isometric actions qd}

\begin{proof}
Consider first a minimal action $\Gamma\acts X$.  Let $G$ be the closure of $\Gamma\subset \text{Isom}(X)$.  Then picking any $x\in X$ and taking the orbit map $G\to Gx$ gives a continuous, equivariant map $h: G\to X$.  Then the map which sends $f\in C(X)$ to $f\circ h$ is an equivariant homomorphism $C(X)\to C(G)$, and minimality and the definition of $h$ imply this map preserves the supremum norm.  The previous theorem now shows $\Gamma\acts X$ is quasidiagonal with completely positive maps as in the previous theorem, as any such map can be pulled back to such a map of $C(X)$.  

If $\Gamma\acts X$ is not minimal, for any fixed $\delta>0$ we can find a finite collection $\{O_i\}$ of orbits so that the union of all their closures, $Y$, is $\delta$-dense in $X$.  Then the map $\psi:C(X)\to C(Y)$ given by restriction is an equivariant, unital homomorphism, and for any finite subset $F\subset C(X)$, $\delta$ can be chosen small enough so that $\psi$ decreases norms by less than any $\epsilon>0$.  Composing with $\psi: C(Y)\to \oplus_i C(O_i)$ and using the previous paragraph now finishes the proof.  \end{proof}

\normalfont

The theorem above essentially implies that crossed products by isometric actions have finite dimensional approximations as nice as those of the acting group's $C^*$-algebra.  More precisely, we have the corollary below, which is an immediate consequence of this theorem and \cite[3.4 and 3.5]{kerr2011residually}.  

\corollary{Suppose $\Gamma\acts X$ is an isometric action by a countable discrete group.  Then $C(X)\rtimes_r \Gamma$ is quasidiagonal or MF exactly when $C_r^*(\Gamma)$ is.  } \qed

\normalfont
This corollary allows us to extend results about the MF property for group $C^*$-algebras to the case of crossed products.  For example, we get the following by using \cite[3.9]{MF_algebras}:

\corollary{Let $\Gamma_i$, $i\in I$ be a countable collection of countable discrete Abelian groups with a common subgroup $H$ and suppose $\star_H \Gamma_i\acts X$ is isometric.  Then $C(X)\rtimes_r \Gamma$ is MF.  }\qed

\bibliography{MyBibliography2.bib}
\bibliographystyle{plain}

\end{document}